\documentclass{amsart}

\usepackage{amsmath}
\usepackage{amscd}
\usepackage{amssymb} 

\newcommand{\cal}{\mathcal}
\newcommand{\bk}{{\bf k}}
\newcommand{\bs}{{\bf s}}
\newcommand{\fg}{{\frak g}}

\newcommand{\wh}{\wedge_h}
\newcommand{\ww}{\wedge_{w}}

\DeclareMathOperator{\End}{End}
\DeclareMathOperator{\Img}{Im}

\DeclareMathOperator{\Ker}{Ker}

\newtheorem{theorem}{Theorem}[section]
\newtheorem{theorem/definition}{Theorem/Definition}[section]
\newtheorem{Theorem}{Theorem}

\newtheorem{Lemma}{Lemma}

\begin{document}

\title
{On Quantum de Rham Cohomology Theory}
\author{Huai-Dong Cao \& Jian Zhou}
\address{Department of Mathematics\\
Texas A\&M University\\
College Station, TX 77843}
\email{cao@math.tamu.edu,  zhou@math.tamu.edu}
\begin{abstract}
We define quantum exterior product $\wedge_h$ 
and quantum exterior differential $d_h$
on  Poisson manifolds (of which symplectic manifolds are an important class of
examples). 
Quantum de Rham cohomology, which is a deformation quantization of de Rham cohomology,
is defined as the cohomology of $d_h$.
We also define quantum Dolbeault cohomology. 
A version of quantum integral on symplectic manifolds is considered and the 
correspoding quantum Stokes theorem is proved. We also derive quantum hard 
Lefschetz theorem.
By replacing $d$ by $d_h$ and 
$\wedge$ by $\wh$ in the usual definitions, 
we define many quantum analogues of
important objects in differential geometry, e.g. 
quantum curvature. 
The quantum characteristic classes are then studied 
along the lines of classical Chern-Weil theory.
Quantum equivariant de Rham cohomology is defined 
in the similar fashion.
\end{abstract}
\maketitle
\date{}
\footnotetext[1]{1991 {\em Mathematics Subject Classification}: }

\footnotetext[2]{Authors' research was supported in part by NSF }

Recently, quantum cohomology has attracted a lot of attentions.
See Tian \cite{Tia} for a survey on this topic, and the introduction of Li-Tian 
\cite{Li-Tia} for more recent development.
In this note we announce a new construction of a deformation 
of the ring  structure on the de Rham cohomology of a closed symplectic
manifold. 
Our construction  is of algebraic nature, 
following the definition of  the de Rham cohomology.
In fact, on any Poisson manifold, we will define a quantum
exterior product $\wh$ of exterior forms, 
and quantum exterior differential $d_h$, such that $d_h ^ 2 = 0$,
and $d_h$  is a derivation for $\wh$. 
We define {\em quantum de Rham cohomology} as the cohomology of $d_h$.
By replacing $d$ by $d_h$ and 
$\wedge$ by $\wh$ in the usual definitions, 
we define many quantum analogues of
important objects, such as quantum curvature and quantum characteristic classes, 
in differential geometry.

Our construction  has the following features which are not
shared by the quantum cohomology:
\begin{enumerate}
\item The proof of the associativity of the quantum exterior product is 
elementary.
\item It is routine to define quantum Dolbeault cohomology.
\item It is routine to define quantum characteristic class.
\item It is routine to define quantum equivariant de Rham cohomology.
\item The computations for homogeneous examples are elementary.
\end{enumerate}

The definition of our quantum exterior product 
is motivated by Moyal-Weyl multiplication and Clifford multiplication.
For any finite dimensional vector space $V$, let  $\{ e_1, \cdots, e_m \}$ be a basis 
of $V$ and $\{ e^1, \cdots, e^m \}$ be the dual basis.
Assume that $w = w^{ij} e_i \otimes e_j \in V \otimes V$.
Then $w$ defines a multiplication $\wedge_w$ on $\Lambda(V^*)$,
and a multiplication $*_w$ on $S(V^*)$,
such that $e^i \wedge_w e^j =   e^i \wedge e^j + w^{ij}$,
$e^i *_w e^j = e^i \odot  e^j + w^{ij}$.
If $w \in S^2(V)$, then $\wedge_w$ is the Clifford multiplication.
If $w \in \Lambda^2(V)$ is nondegenerate, $*_w$ is 
the Moyal-Weyl multiplication.
If $w \in \Lambda^2(V)$, then $\wedge_w$ is what we call the  {\em quantum 
exterior product} (or a quantum Clifford multiplication).
It is elementary to show that this multiplication is associative. 
We use it to obtain a  quantum calculus on  any Poisson manifold.
We will present the main results below, and the proofs will appear elsewhere.

\vspace{.2in}

\noindent {\bf 1. Quantum exterior algebra}. 
Let $V$ be a finite dimensional vector space over a 
field $\bk$ of characteristic zero,
and $\Lambda(V^*)$ be the exterior algebra generated by the dual vector space
$V^*$. 
For any $v \in V$ and $\alpha \in \Lambda^k(V^*)$, denote
\begin{eqnarray*}
(v \vdash \alpha ) (v_1, \cdots, v_{k-1}) & = &
	 \alpha (v, v_1, \cdots, v_{k-1}), \\
(\alpha \dashv v) (v_1, \cdots, v_{k-1}) & = & 
	\alpha (v_1, \cdots, v_{k-1}, v),
\end{eqnarray*}
for $v_1, \cdots, v_{k-1} \in V$.
Let $\Lambda_h(V^*) = \Lambda(V^*)[h] = \Lambda(V^*) \otimes_{\bk} \bk[h]$.
For any  $w \in \Lambda^2(V)$, with $w = w^{ij}e_i \wedge e_j$ with respect to a basis
$\{ e_1, \cdots, e_m \}$ of $V$, we define the {\em quantum exterior product} 
$\wedge_{h,w}: \Lambda(V^*) \otimes \Lambda(V^*) \rightarrow 
	\Lambda(V^*)[h]$ by
\begin{eqnarray*}
&   & \alpha \wedge_{h, w} \beta \\
& = & \sum_{n \geq 0} 
	\frac{h^n}{n!} w^{i_1j_1} \cdots w^{i_nj_n}
	(\alpha \dashv e_{i_1} \dashv \cdots \dashv e_{i_n})
	\wedge (e_{j_n} \vdash \cdots \vdash e_{j_1} \vdash \beta), 
\end{eqnarray*}	
for $\alpha, \beta \in \Lambda(V^*)$.
This definition is evidently independent of the choice of the basis.
We extend $\wedge_{h, w}$ as $\bk[h]$-module map to 
$\Lambda_h(V^*) \otimes_{\bk[h]} \Lambda_h(V^*)$.
When there is no confusion about $w$, we will simply write
$\alpha \wh \beta$ for $\alpha \wedge_{h, w} \beta$.
We are interested in $\alpha \wedge_w \beta$, which is just 
$\alpha \wedge_{1, w} \beta$.
We assign $h$ the degree $2$. Then $\Lambda_h(V^*)$ 
has a natural $\Bbb Z$-grading.
Denote by $\Lambda_h^{[n]}(V^*)$ the subspace of homogeneous elements
of degree $n$, then it is clear that
$$\Lambda_h^{[m]}(V^*) \wedge_h \Lambda_h^{[n]}(V^*) \subset
	\Lambda_h^{[m+n]}(V^*).
$$

\begin{Theorem} (Cao-Zhou \cite{Cao-Zho}) \label{Thm:algebraic}
The quantum exterior product satisfies the following properties
\begin{eqnarray}
\text{Supercommutativity} 
& & \alpha \wh \beta = (-1)^{|\alpha||\beta|} \beta \wh \alpha, 
	\label{com} \\
\text{Associativity} 
& & (\alpha \wh \beta) \wh \gamma = \alpha \wh (\beta \wh \gamma), 
	\label{ass}
\end{eqnarray}
for  all $\alpha, \beta, \gamma \in \Lambda_h(V^*)$.
Therefore, $(\Lambda_h(V^*), \lambda_h)$ is a deformation quantization
of the exterior algebra $(\Lambda(V^*), \wedge)$.
\end{Theorem}

The proof of $(\ref{com})$ is trivial. The proof of $(\ref{ass})$ is of 
elementary nature but non-trivial.  It is proved first by brute force in the case of $\deg (\alpha) = 1$,
then by induction on $\deg (\alpha)$ (see \cite{Cao-Zho} for details).

We can  also extend $\wh$ to $\Lambda_{h, h^{-1}}(V^*) 
= \Lambda(V^*)[h, h^{-1}] = \Lambda(V^*) \otimes_{\bk} \bk[h, h^{-1}]$. 

An algebra $A$ with unit $e \in A$ 
over a field $\bk$ is called a {\em $\bk$-Fr\"{o}benius algebra} 
if there is a nondegenerate symmetric 
$\bk$-bilinear function
$<\cdot, \cdot>: A \times A \rightarrow \bk$,
such that
$$<\alpha \beta, \gamma> = <\alpha, \beta \gamma>,
$$
for all $\alpha, \beta, \gamma \in A$.
There is a simple way to construct a structure of Fr\"{o}benius 
algebra on any $\bk$-algebra $A$ with unit. 
Let $\phi: A \rightarrow \bk$ be a nonzero $\bk$-functional 
on $A$.
Set $<\alpha, \beta>_{\phi} = \phi(\alpha\beta)$ 
for $\alpha, \beta \in A$.
If it is nondegerate, then $(A, <\cdot, \cdot>_{\phi})$ is
a $\bk$-Fr\"{o}benius algebra.
Conversely, given any Fr\"{o}benius algebra 
$(A, <\cdot, \cdot>)$,
let $\phi(\alpha) = <\alpha, e>$, for $\alpha \in A$,
then $<\cdot, \cdot> = <\cdot, \cdot>_{\phi}$.
Now on $\Lambda(V^*)$, consider a Berezin integral 
$\int: \Lambda(V^*) \rightarrow \bk$ ( i.e., a $\bk$-linear
functional which is only nonzero on $\Lambda^{top}(V^*)$).
Then it is clear that $<\alpha, \beta> = \int \alpha \wedge_w \beta$
defines a structure of Fr\"{o}benius algebra on $(\Lambda(V^*), \wedge_w)$.

\vspace{.2in}

\noindent {\bf 2. Quantum de Rham complex}. 
Let $(P, w)$ be a  Poisson manifold, with bi-vector field $w$, whose Schouten-Nijenhuis bracket
vanishes. 
(See Vaisman \cite{Vai} for definitions.) 
The fiberwise quantum exterior multiplication defines 
\begin{eqnarray*}
& \wh: \Omega_h(M) \otimes \Omega_h(M) \rightarrow \Omega_h(M), \\
& \wh: \Omega_{h, h^{-1}}(M) \otimes \Omega_{h, h^{-1}}(M) 
	\rightarrow \Omega_{h, h^{-1}}(M),
\end{eqnarray*}
where $\Omega_{h}(V^*) 
= \Omega(V^*)[h] = \Omega(V^*) \otimes_{\bk} \bk[h]$
and $\Omega_{h, h^{-1}}(V^*) 
= \Omega(V^*)[h, h^{-1}] = \Omega(V^*) \otimes_{\bk} \bk[h, h^{-1}]$.
Koszul \cite{Kos} defined an operator 
$\delta: \Omega^{k}(M) \rightarrow \Omega^{k-1}(M)$ by
$$\delta \alpha = w \vdash d \alpha - d (w \vdash \alpha),$$
for $\alpha \in \Omega^k(M)$.
He also showed that $\delta^2 = 0$, 
$d\delta + \delta d = 0$.
We define the quantum exterior differential $d_h = d - h \delta: \Omega(M)[h] 
\rightarrow
\Omega(M)[h]$ , and similarly on $\Omega(M)[h,h{^-1}]$.
Then it is easy to see that $d_h^2 = 0$.
One of the technical results in Cao-Zhou \cite{Cao-Zho} is the following

\begin{Theorem} \label{Thm:main}
For any Poisson manifold $(M, w)$, 
$d_h$ satisfies 
\begin{eqnarray} \label{derivation2}
d_h(\alpha \wh \beta) = 
	(d_h \alpha) \wh \beta + 
	(-1)^{|\alpha|}\alpha \wh (d_h \beta), 
\end{eqnarray}
for $\alpha, \beta$ in $\Omega(M)[h]$,  or 
 $\alpha, \beta$ in $\Omega(M)[h, h^{-1}]$. 
\end{Theorem}

This is proved first for $\deg (\alpha) = 1$ by brute force,
then by induction  on $\deg (\alpha)$.
For regular Poisson manifolds (e.g., symplectic manifolds),
there is an easier proof.
On such Poisson manifolds, there always exists a torsion-less connection
$\nabla$ which preserves $w$.
Then for any local frame $\{ e^1, \cdots, e^n \}$, 
and $\alpha \in \Omega(M)$, we have 
$$d_h \alpha = e^i \wh \nabla_{e_i} \alpha.$$
This is the analogue of a similar expression for $d + d^*$ in
Riemannian geometry (Lawson-Michelsohn \cite{Law-Mic}, Lemma II.5.13).
Use the analogue of normal coordinates, 
the proof of Theorem \ref{Thm:main} reduces to 
the associativity of the quantum exterior multiplication.

\vspace{.2in}

\noindent {\bf 3. Quantum de Rham cohomology}.
For any Poisson manifold $(M, w)$, 
the {\em (polynomial) quantum deRham cohomology} is defined
by 
$$Q_hH_{dR}^*(M) = \Ker d_h / \Img d_h,$$
for the quantum exterior differential $d_h: \Omega(M)[h] \rightarrow \Omega(M)[h]$.
The {\em Laurent quantum de Rham cohomology} is
$$Q_{h, h^{-1}}H_{dR}^*(M) = \Ker d_h / \Img d_h,$$
for  $d_h: \Omega(M)[h, h^{-1}] \rightarrow 
\Omega(M)[h, h^{-1}]$.
As a consequence of Theorem \ref{Thm:algebraic} and Theorem \ref{Thm:main}, 
we have
\begin{Theorem} The quantum de Rham cohomology 
$Q_hH_{dR}^*(M)$  of a Poison manifold  
has the following properties:
\begin{eqnarray*}
\alpha \wh \beta & = & (-1)^{|\alpha||\beta|}\beta \wh \alpha, \\
(\alpha  \wh \beta) \wh \gamma & = & \alpha \wh (\beta \wh \gamma),
\end{eqnarray*}
for $\alpha, \beta, \gamma \in Q_hH^*_{dR}(M)$.
Similar results hold for Laurent quantum deRham cohomology.
\end{Theorem}

The complex $(\Omega(M)[h], d_h)$ can be regarded as a double complex
$(C^{p, q}, -h\delta, d)$, where
$C^{p,q} = h^p\Omega^{q-p}(M)$, $p \geq 0$.
This is the analogue of Brylinski's double complex ${\cal C}_{..}(M)$ 
(\cite{Bry}, $\S 1.3$).
By the standard theory for double complex
(Bott-Tu \cite{Bot-Tu}, $\S 14$), 
there are two spectral sequences $E$ and $E'$ abutting to 
$H^*(\Omega[h], d_h) = Q_hH^*_{dR}(M)$,
with $E_1^{p, q} = h^pH^q(C^{p, *}, d) = h^pH^{q-p}_{dR}(M)$,
$(E_1')^{p, q} = h^pH^*(C^{*, q}, \delta) = h^pPH_{q-p}(M)$, $p \geq 0$. 
Since nontrivial $E_1^{p, q}$ all have $p+ q$ even,
and the differential $d_r$ change the parity of $p+q$, 
it is routine to prove the following

\begin{Theorem} \label{Thm:spectral1}
For a Poisson manifold with odd Betti numbers all vanishing,
the spectral sequence $E$ degenerate at $E_1$, i.e. $d_r = 0$ for
all $r \geq 0$, hence
$Q_hH^*_{dR}(M)$ is a deformation quantization of $H^*_{dR}(M)$.
\end{Theorem}

Brylinski \cite{Bry} proved that on closed K\"{a}hler manifold $(M, \omega)$,
every de Rham cohomology class has a representative $\alpha$ such that
$d \alpha = 0$, $\delta \alpha = 0$.
This implies that 

\begin{Theorem} 
For a closed K\"{a}hler manifold $M$,
the spectral sequence $E$ degenerate at $E_1$, i.e. $d_r = 0$ for
all $r \geq 0$, hence
$Q_hH^*_{dR}(M)$ is a deformation quantization of $H^*_{dR}(M)$.
\end{Theorem}

Similarly, we regard $(\Omega(M)[h, h^{-1}], d - h \delta)$ 
 as a double complex 
$(\widetilde{C}^{p, q}, -h\delta, d)$,
where $\widetilde{C}^{p, q} = h^p\Omega^{q-p}(M)$, 
$p, q \in {\Bbb Z}$. 
This is essentially Brylinski's double complex 
${\cal C}_{..}^{per}$, but with a different bi-grading.
We get two spectral sequences $\tilde{E}$ and  $\tilde{E}'$
abutting to $Q_{h, h^{-1}}H^*_{dR}(M)$, 
with $\tilde{E}_1^{p, q} = h^pH^{q-p}_{dR}(M)$, 
$(\tilde{E}_1')^{p, q} = h^pH^*(C^{*, q}, \delta) = h^pPH_{q-p}(M)$,
$p, q \in {\Bbb Z}$. 
It is clear that analogue of Theorem \ref{Thm:spectral1}
holds for $\tilde{E}$.
On the other hand, by a method of Brylinski \cite{Bry}, 
one can prove the following:

\begin{Theorem} \label{thm:degeneracy}
For any compact symplectic manifold without boundary,
the spectral sequences $\tilde{E}$ and $\tilde{E}'$ degenerate 
at $\tilde{E}_1$ and $\tilde{E}_1'$ respectively.
Hence
$Q_{h, h^{-1}}H^*_{dR}(M)$ is a Laurent 
deformation quantization of $H^*_{dR}(M)$.
\end{Theorem}

Fixing an isomorphism $H^{2n}_{dR}(M) \cong {\Bbb R}$ then defines
a structure of Fr\"{o}benius algebra on $(H^*_{dR}(M), \wedge_w)$.

\vspace{.2in}

\noindent {\bf 4. Quantum Hard Lefschetz Theorem}.
For a closed symplectic manifold, the analogue of Hard Lefschetz
Theorem \cite{Gri-Har} holds for $Q_{h, h^{-1}}H^*_{dR}(M)$.
We begin with a $2n$-dimensional symplectic vector space $(V, \omega)$.
Brylinski \cite{Bry} defined a symplectic star operator 
$*: \Lambda^k(V^*) \rightarrow \Lambda^{2n-k}(V^*)$.
We can extend it to $\Lambda_{h, h^{-1}}$ by setting $*h = h^{-1}$,
and $*h^{-1} = h$.
Define $L_h: \Lambda_{h, h^{-1}}(V^*) \rightarrow \Lambda_{h, h^{-1}}(V^*)$
by $L_h(\alpha) = \omega \wh \alpha$.
Define $L_h^* = -*L*$, and $A_h: \Lambda_{h, h^{-1}}(V^*) \rightarrow 
	\Lambda_{h, h^{-1}}(V^*)$
by $A_h(\alpha) = (n - k) \alpha$, 
for $\alpha \in \Lambda^{[k]}_{h, h^{-1}}(V^*)$.
Then we have

\begin{Lemma} 
The following identities hold:
\begin{eqnarray*}
[L_h, L_h^*] = 0, & [L_h, A_h] = 2 L_h, & [L^*_h, A_h] = - 2 L^*_h.
\end{eqnarray*}
Furthermore, if we regard multiplications by $h$ and $h^{-1}$
as operators, then we have 
\begin{eqnarray*}
[h, h^{-1}] = 0, & [L_h, h^{\pm 1}] = [L_h^*, h^{\pm 1}] = 0, &
[A_h, h^{\pm}] = \pm 2 h^{\pm 1}.
\end{eqnarray*}
\end{Lemma}

Thus we cannot use the representation theory of $sl(2, {\Bbb C})$
as in the classical theory.
Notice that multiplication by $h$ is an isomorphism, whose inverse is
multiplication by $h^{-1}$.
Let $M_h = h^{-1}L_h$, $M_h^* = hL_h^*$, then $M_h, M_h^*, A_h$
form an abelian algebra. 
Notice that it now suffices to find the eigenvalues of $M_h$ 
on $\Lambda^{[0]}_{h, h^{-1}}(V^*)$ and $\Lambda^{[1]}_{h, h^{-1}}(V^*)$.
In Cao-Zhou \cite{Cao-Zho},
this is done by induction on the dimension of $V$ and the following 

\begin{Lemma}
 Let $\{ M_n \}$ be a sequence of square 
matrices with coefficient in $\bk$
obtained in the following way:
$$
M_{n+1} = \left( \begin{array}{cc} 
M_n & -I \\ I & M_n + 2I \end{array} \right)
$$
for $n \geq 1$. 
where $I$ is the identity matrix of the same size as $M_n$.
For any $\lambda \in \bk$, and $n \geq 1$, we have
$$\det (M_{n+1} + \lambda I)  
=  \det [M_n + (\lambda +1)I]^2.$$
Therefore, the eigenvalues of $M_{n+1}$ can be 
obtained by adding $1$ to that of $M_n$,
with the multiplicities doubled.
\end{Lemma}

As a consequence, the eigenvalues of $M_h$ on $\Lambda_{h, h^{-1}}(V^*)$
are $n$ and $n \pm \sqrt{5}/2$, when $\dim (V) = 2n$.
For $M_h^*$, they are $-n$ and $-n \pm \sqrt{5}/2$. Therefore, we have

\begin{Theorem} 
For a symplectic vector space $V$, the operators $L_h$ and
$L_h^*$ are  isomorphisms.
Furthermore, $\Lambda_{h, h^{-1}}(V^*)$ decomposes 
into one dimensional eigenspaces of $h^{-1}L_h$ 
(or $hL_h^*$) with nonzero eigenvalues.
\end{Theorem}

\begin{Lemma} 
On a symplectic manifold $(M, \omega)$, we have
\begin{eqnarray*}
[L_h, d_h] = 0, & [L_h^*, d_h] = 0, & [A_h, d_h] = -d_h.
\end{eqnarray*}
\end{Lemma}

As a consequence, we get:

\begin{Theorem} (Quantum Hard Lefschetz Theorem) 
For any symplectic manifold $(M^{2n},  \omega)$,
its Laurent quantum deRham cohomology 
$Q_{h, h^{-1}}H^*_{dR}(M)$ decomposes into 
one-dimensional eigenspaces of  the operator $h^{-1}L_h$
(or $hL_h^*$) with nonzero eigenvalues. 
In particular, $L_h$ and $L_h^*$ are isomorphisms.
\end{Theorem}

\vspace{.2in}

\noindent {\bf 5. Complexified quantum exterior algebra}.
We also consider real vector space $V$ with an
almost complex structure $J\in End (V)$
 such that $J^2 = - Id$.
There is an induced linear transformation 
$\Lambda^2J: \Lambda^2(V) \rightarrow \Lambda^2(V)$.
For any bi-vector $w \in \Lambda^2(V)$,
$J$ is said to preserve $w$ if 
$\Lambda^2J (w) = w$.
Given any bi-vector $w$ which is preserved by $J$, 
we can define the quantum exterior product on 
$\Lambda_h(V^*)$ as in the last section.
Now if we tensor everything by $\Bbb C$,
we get a complex algebra ${\Bbb C}\Lambda_h(V^*)$, which is
a deformation quantization of ${\Bbb C}\Lambda(V^*) :=
\Lambda(V^*) \otimes_{\Bbb R} {\Bbb C} 
= \Lambda_{\Bbb C}(V^* \otimes_{\Bbb R} {\Bbb C})$.  As in complex geometry, 
we can exploit a natural decomposition as follows. $J$ can be uniquely extended to a 
complex linear endomorphism, denoted also by $J$, of  ${\Bbb C}V$ also satisfying 
$J^2 = - Id$. There is a natural identification of complex vector
spaces  ${\Bbb C}V \cong V^{1, 0} \oplus V^{0, 1}$ , where $V^{1, 0}$ and $V^{0, 1}$
are eigenspaces of $J$ with eigenvalues $\sqrt {-1}$ and $-\sqrt {-1}$ respectively.
As a consequence, there are decompositions
\begin{eqnarray*}
{\Bbb C}\Lambda(V) & = & \oplus_{p, q} \Lambda^{p, q}(V), \\
{\Bbb C}\Lambda(V^*) & = & \oplus_{p, q} \Lambda^{p, q}(V^*),
\end{eqnarray*}
where $\Lambda^{p, q}(V) \cong \Lambda_{\Bbb C}^p(V^{1, 0})
\otimes_{\Bbb C} \Lambda^q_{\Bbb C}(V^{0, 1})$,
and $\Lambda^{p, q}(V^*) \cong \Lambda_{\Bbb C}^p((V^{1, 0})^*)
\otimes_{\Bbb C} \Lambda^q_{\Bbb C}((V^{0, 1})^*)$.
We give ${\Bbb C}\Lambda_h(V^*)$ the following 
${\Bbb Z} \times {\Bbb Z}$-bigrading: 
elements in $\Lambda^{p, q}(V^*)$
has bi-degree $(p, q)$,  and $h$ has bi-degree $(1, 1)$.
Since $w$ is preserved by $J$, it belongs to $\Lambda^{1, 1}(V)$
after complexification.
Denote by $\Lambda_h^{[p, q]}(V^*)$ the space of homogeneous
elements of bi-degree $(p, q)$.
It is then straightforward to see that
$$\Lambda_h^{[p, q]}(V^*) \wh \Lambda_h^{[r, s]}(V^*)
	\subset \Lambda_h^{[p+r, q+t]}(V^*).$$

Now let  $\omega$ be a symplectic 
form on $V$ which is compatible with an almost complex structure
on $V$. Namely, rank of $\omega$ is $2n = \dim (V)$, 
$w(J\cdot, J\cdot) = \omega(\cdot, \cdot)$,
and $g(\cdot, \cdot) := \omega(\cdot, J\cdot)$ is 
a positive definite element of $S^2(V^*)$.
Then $\omega$ induces a natural Hermitian metric $H$ on 
${\Bbb C}\Lambda(V^*)$, such that
$$ H(\alpha \ww \beta, \gamma)  
= H(\alpha, \beta \ww \gamma)$$
for any $\alpha, \beta, \gamma \in {\Bbb C}\Lambda(V^*)$.
Here $w \in \Lambda^2(V)$ is obtained from $\omega$ by 
``raising the indices" (For details, see Cao-Zhou \cite{Cao-Zho}, $\S 1.5$). 
This shows that the algebra $({\Bbb C}\Lambda(V^*), \wedge_w)$ 
has a structure of Hermitian Fr\"{a}benius algebra.

\vspace{.2in}

\noindent {\bf 6. Quantum Dolbeault cohomology}.

On a complex manifold $(M, J)$ with a Poisson structure $w$, 
such that $J$ preserves $w$,
define $\delta^{-1, 0}: \Omega^{p, q}(M) 
	\rightarrow \Omega^{p-1, q}(M)$
and $\delta^{0, -1}(M): \Omega^{p, q}(M) 
	\rightarrow \Omega^{p, q-1}(M)$ by
\begin{eqnarray*}
\delta^{0, -1} \alpha  & = & w\vdash (\partial \alpha) 
	- \partial (w \vdash \alpha), \\
\delta^{-1, 0} \alpha  & = & 
	w\vdash (\overline{\partial} \alpha) 
	- \overline{\partial} (w \vdash \alpha),
\end{eqnarray*}
for $\alpha \in \Omega^{p, q}(M)$.
Set $\partial_h = \partial - h \delta^{0, -1}$,
and $\overline{\partial}_h = \overline{\partial} - h \delta^{-1, 0}$.
Then $d_h = \partial_h + \overline{\partial}_h$.
Now $0 = d_h^2 = \partial_h^2 + 
	(\partial_h\overline{\partial}_h 
	+ \overline{\partial}_h\partial_h)
	+ \overline{\partial}_h^2$,
since they have bi-degrees $(2, 0)$, 
$(1, 1)$ and $(0, 2)$ respectively.
Hence, we have
\begin{eqnarray*} \label{double}
\partial_h^2 = 0, & 
\partial_h\overline{\partial}_h 
	+ \overline{\partial}_h\partial_h = 0, &
\overline{\partial}_h^2 = 0.
\end{eqnarray*}
We then define
\begin{eqnarray*}
Q_hH^{p, *}(M) & = & H(\Omega_h^{[p, *]}(M), \overline{\partial}_h), \\
Q_{h, h^{-1}}H^{p, *}(M) & = & 
	H(\Omega_{h, h^{-1}}^{[p, *]}(M), \overline{\partial}_h).
\end{eqnarray*}
They wil be called {\em quantum Dolbeault cohomology} and
{\em Laurent quantum Dolbeault cohomology} respectively.
Several relevant spectral sequences and their degeneracy are considered
in Cao-Zhou \cite{Cao-Zho}.

\vspace{.2in}

\noindent {\bf 7. Quantum integral and quantum Stokes Theorem}

Let $(M, \omega)$ be a closed $2n$-dimensional symplectic manifold.
Define an integral $\int_h: \Omega_h(M) \rightarrow {\Bbb R}[h]$
as follows.
For any $\alpha \in \Omega^{j}(M)$,
if $j$ is odd, set $\int_h \alpha = 0$;
if $j = 2n - 2k $ for some integer $k$,
set 
$$\int_h \alpha = \int_M \alpha \wedge \frac{\omega^k}{k!}.$$
Extend $\int_h$ to $\Omega_h(M)$ as a ${\Bbb R}[h]$-module map.
We call $\int_h$ the {\em quantum integral}.
Then we have

\begin{theorem} (Quantum Staokes Theorem) 
For any $\alpha \in \Omega^{j}(M)$, we have $\int_h d \alpha = 0$, 
        $\int_h h\delta \alpha = 0$, and therefore
$$\int_h d_h \alpha = 0.$$
\end{theorem}

\vspace{.2in}

\noindent {\bf 8. Quantum Chern-Weil theory}.

Given a real or complex vector bundle $E \rightarrow M$  over a 
Poisson manifold $M$, and a connection $\nabla^E$ on it,
we define the quantum covariant derivative
$$d_h^{\nabla_E}: \Omega_h^*(E) \rightarrow \Omega_h^*(E)$$
 as follows.
Let $\bs$ be a local frame of $E$ and $\theta$ be the 
connection $1$-form in this frame. Then we have: 
$\nabla \bs = \bs \otimes \theta$, i.e.,
$$\nabla \bs_j = \sum_{k=1}^{n} \bs_k \otimes \theta_j^k.$$
For $\alpha =\bs \otimes \phi$, where $\phi$ is some vector valued form,
we define
$$d_h^{\nabla_E} \alpha = \bs \otimes (\theta \wh \phi + d_h \phi) 
	= \sum  \bs_k \otimes ( \theta_j^k \wh \phi^j + d_h \phi^k).$$
It is straightforward to check 
that the definition of $d_h^{\nabla_E}$ is independent of 
the choice of the local frames (Cao-Zhou \cite{Cao-Zho}, Lemma 6.1).
Furthermore, there is an element $R^E_h \in
\Omega^2_h(\End(E))$, such that
for each $k \geq 0$, $(d_h^{\nabla_E})^2$ on $\Omega^k_h(M)$ is
given by $(d_h^{\nabla_E})^2 \Phi = \Phi \wh R_h^E$, for
any $\Phi \in \Omega^*_h(E)$.
$R^E_h$ is called the quantum curvature of $\nabla^E$.
In a local frame, $R^E_h$ is given by
$$F_h = d_h \theta + \theta \wh \theta,$$
where $\theta$ is the connection $1$-form matrix in the local frame. 

If $p$ is a polynomial on the space of $n\times n$-matrices, 
such that $p(G^{-1}AG) = p(A)$, 
for any $n \times n$ matrix $A$, and
invertible $n \times n$-matrix $G$,
then $p(F_h)$ for different frames patch up to 
a well-defined element $p(R^E) \in \Omega^*(M)[h]$.
Similar to the ordinary Chern-Weil theory, 
it is easy to see that $d_h p(R^E) = 0$.
So it defines a class in $Q_hH^*_{dR}(M)$.
The usual  construction of transgression operator  
carries over to show that
this class is independent of the choice of the connection $\nabla^E$.
In this way,  one can define quantum Chern classes, 
quantum Euler class, etc.
We will call them quantum characteristic classes.
It is clear that we can repeat the same story in Laurent case.
\vspace{.2in}

\noindent {\bf 9. Quantum equivariant de Rham cohomology}.
Let $(M, w)$ be a Poisson manifold, which 
admits an action by a compact connected Lie group $G$,
such that the $G$-action preserves the Poisson bi-vector field $w$.
Let $\fg$ be the Lie algebra of $G$, $\{ \xi_a \}$ a basis of $\fg$ and 
$\{ \Theta^a \}$ be the dual basis in $S^1(\fg^*)$. 
Denote by $\iota_a$ the contraction by the vector field 
generated by the one parameter group corresponding to $\xi_a$,
and $L_a$ the Lie derivative by the same vector field.
Imitating the Cartan model for equivariant cohomology, 
we consider the operator $D_{hG} = d_h + \Theta^a \iota_a 
= d - h \delta + \Theta^a \iota_a$
acting on $(S(\fg^*) \otimes \Omega(M))^G[h]$.  
It is well-known that $d + \Theta^a \iota_a$ maps 
$(S(\fg^*) \otimes \Omega(M))^G$ to itself.
Since the $G$-action preserves $w$, it is easy to check that
$\delta$ also preserves $(S(\fg^*) \otimes \Omega(M))^G$.
Therefore, $D_{hG}$ is an operator from 
$(S(\fg^*) \otimes \Omega(M))^G[h]$
to itself. 
Now on $(S(\fg^*) \otimes \Omega(M))^G[h]$, we have 
\begin{eqnarray*}
D_{hG}^2 
& = & d_h^2 + (\Theta^a\iota_a)^2 + \Theta^a (d \iota_a + \iota_a d)
	- h \Theta^a (\delta \iota_a + \iota_a \delta) \\
& = &	- h \Theta^a (\delta \iota_a + \iota_a \delta).
\end{eqnarray*}
Since $\delta = \iota_{w} d - d \iota_{w}$, we have
\begin{eqnarray*}
\delta \iota_a + \iota_a \delta 
& = & \iota_{w}d \iota_a - d \iota_{w} \iota_a 
	+ \iota_a \iota_{w} d - \iota_a d \iota_{w} \\
& = & \iota_{w}d \iota_a - d \iota_a \iota_{w}
	+ \iota_{w} \iota_a d - \iota_a d \iota_{w} \\
& = & \iota_{w}L_a - L_a \iota_{w} = -\iota_{L_aw} = 0.
\end{eqnarray*}
Hence, $D_{hG}^2 =0$.
We call the cohomology 
$$Q_hH^*_{G}(M) := H^*((S(\fg^*) \otimes \Omega(M))^G[h], D_{hG})$$
the quantum equivariant de Rham cohomology. 
Similar definitions can be made using Laurent deformation. 
We will study quantum equivariant de Rham cohomology in a forthcoming
paper.

\end{document}